\documentclass[11pt]{article}

\title
{Finding and using expanders in locally sparse graphs}
\author{Michael Krivelevich
\thanks{School of Mathematical Sciences, Raymond and Beverly
Sackler Faculty of Exact Sciences, Tel Aviv University, Tel Aviv,
6997801, Israel. Email: krivelev@post.tau.ac.il. Research supported in
part by USA-Israel BSF grant 2014361.}
}

\usepackage{amsmath,amsthm, amssymb,latexsym, amsfonts}

\begin{document}
\bibliographystyle{plain}
\maketitle
\newtheorem{thm}{Theorem}
\newtheorem{propos}{Proposition}
\newtheorem{defin}{Definition}
\newtheorem{lemma}{Lemma}[section]
\newtheorem{corol}{Corollary}
\newtheorem{thmtool}{Theorem}[section]
\newtheorem{corollary}[thmtool]{Corollary}
\newtheorem{lem}[thmtool]{Lemma}
\newtheorem{prop}[thmtool]{Proposition}
\newtheorem{clm}[thmtool]{Claim}
\newtheorem{conjecture}{Conjecture}
\newtheorem{problem}{Problem}
\newcommand{\Proof}{\noindent{\bf Proof.}\ \ }
\newcommand{\Remarks}{\noindent{\bf Remarks:}\ \ }
\newcommand{\Remark}{\noindent{\bf Remark:}\ \ }
\newcommand{\whp}{{\bf whp}\ }
\newcommand{\prob}{probability}
\newcommand{\rn}{random}
\newcommand{\rv}{random variable}
\newcommand{\hpg}{hypergraph}
\newcommand{\hpgs}{hypergraphs}
\newcommand{\subhpg}{subhypergraph}
\newcommand{\subhpgs}{subhypergraphs}
\newcommand{\bH}{{\bf H}}
\newcommand{\cH}{{\cal H}}
\newcommand{\cT}{{\cal T}}
\newcommand{\cF}{{\cal F}}
\newcommand{\cD}{{\cal D}}
\newcommand{\cC}{{\cal C}}

\begin{abstract}
We show that every locally sparse graph contains a linearly sized expanding subgraph. For constants $c_1>c_2>1$, $0<\alpha<1$, a graph $G$ on $n$ vertices is called a $(c_1,c_2,\alpha)$-graph if it has at least $c_1n$ edges, but every vertex subset $W\subset V(G)$ of size $|W|\le \alpha n$ spans less than $c_2|W|$ edges. We prove that every $(c_1,c_2,\alpha)$-graph with bounded degrees contains an induced expander on linearly many vertices. The proof can be made algorithmic.

We then discuss several applications of our main result to random graphs, to problems about embedding graph minors, and to positional games.
\end{abstract}

\section{Introduction and main result}\label{sec-intro}
The main goal of this paper is to introduce, and then to apply, a simple sufficient condition, guaranteeing the existence of a large expanding subgraph in a given graph.

Expanders (see \cite{HLW06,KS06} for two extensive surveys on the subject -- or \cite{Sar04} for its very concise introduction) have become one of the most central, and also one of the most applicable, notions of modern combinatorics. Given the utmost importance of expanders and their wide applicability, it is only natural to witness a very substantial research aiming to provide sufficient conditions for a graph being an expander. Frequently spectral properties are used to guarantee expansion (see \cite{Alo86} for the cornerstone contribution in this direction).

Obviously not every graph is an expander; moreover, the standard notion of expansion is rather fragile -- adding a single isolated vertex to a strong expander $G$ produces a non-expanding graph $G'$. This motivates us to pursue a different avenue here, namely, to try and find a (preferably simple) sufficient condition guaranteeing the existence of a large expanding subgraph in a given graph $G$. Strictly speaking, this research direction is certainly not new, see, e.g., \cite{KS94,KS96}, or \cite{SS15, Mon15} for recent results. (We will discuss them briefly after having introduced our main result later in the paper.) Another related line of research, quite popular in the computer science community, is to decompose a given graph, perhaps after slight alterations, into expanding subgraphs, see, e.g., \cite{GR99, KVV04, Tre08,ABS15} for results of this type.

Let us lay a formal ground to state our result.

As usual, for a graph $G=(V,E)$ and a vertex set $W\subset V$ we denote by $N_G(W)$ the external neighborhood of $W$ in $G$, i.e.,
$$
N_G(W)=\{v\in V\setminus W: v\mbox{ has a neighbor in }W\}\,.
$$

We can now give the formal definition of an expander we adapt in this paper.
\begin{defin}\label{expander}
Let $G=(V,E)$ be a graph on $n$ vertices, and let $\gamma>0$. The graph $G$ is a {\em $\gamma$-expander} if $|N_G(W)|\ge \gamma |W|$ for every vertex subset $W\subset V$, $|W|\le n/2$.
\end{defin}

This is a fairly commonly used notion of an expander, see, e.g., \cite{Alo86}, or \cite[Chapter 9]{AS}. It does not aim to capture or to reflect the strongest possible level of expansion, but it is strong enough to derive many nice graph properties. Here is a short and somewhat informal list, part of which will be  discussed later in the paper: connectedness of $G$; logarithmic diameter; logarithmic mixing time of a random walk on $G$; all separators in $G$ are linearly large; one can embed large minors in $G$. The restriction $|W|\le n/2$ above is rather arbitrary, it reflects the fact that the size of vertex subsets under consideration should be capped from above to allow them room to expand externally in $G$.

As we have promised, we will present a simple sufficient condition for the existence of a large expanding subgraph in a given graph. This condition is based on what we call local sparseness. Here is a formal definition:

\begin{defin}\label{sparseness}
Let $c_1>c_2>1$, $0<\alpha<1$. A graph $G=(V,E)$ on $n$ vertices is called a {\em $(c_1,c_2,\alpha)$-graph} if
\begin{enumerate}
\item
$\frac{|E|}{|V|}\ge c_1$\,;
\item every vertex subset $W\subset V$ of size $|W|\le \alpha n$ spans less than $c_2|W|$ edges.
\end{enumerate}
\end{defin}
In words, the above condition says that relatively small sets are sizably sparser than the whole graph. It has been used in recent paper \cite{Kri16} of the author.

How natural or common is this condition? As the (very easy) proposition below shows, most sparse graphs are locally sparse.

\begin{propos}\label{prop1}
Let $c_1>c_2>1$ be reals. Define $\alpha=\left(\frac{c_2}{5c_1}\right)^{\frac{c_2}{c_2-1}}$. Let $G$ be a random graph drawn from the probability distribution $G\left(n,\frac{c_1}{n}\right)$. Then \whp every set of $k\le \alpha n$ vertices of $G$ spans fewer than $c_2k$ edges.
\end{propos}

One can also easily cap maximum degree in (a nearly spanning subgraph of) a random graph, as given by the following proposition.

\begin{propos}\label{prop2}
For every $C>0$ and all sufficiently small  $\delta>0$ the following holds. Let $G$ be a random graph drawn from the probability distribution $G\left(n,\frac{C}{n}\right)$. Then \whp every set of $\frac{\delta}{\ln\frac{1}{\delta}} n$ vertices of $G$ touches fewer than $\delta n$ edges.
\end{propos}

Observe that if $G=(V,E)$ satisfies the conclusion of the above proposition, then by deleting $\frac{\delta}{\ln\frac{1}{\delta}} n$ vertices of highest degrees in $G$, one obtains a spanning subgraph $G'=(V',E')$ on $|V'|=\left(1-\frac{\delta}{\ln\frac{1}{\delta}}\right) n$ vertices and with $|E'|\ge |E|-\delta n$ edges, in which all degrees are at most $2\ln(1/\delta)$. (Otherwise, all deleted vertices are of degree at least $2\ln(1/\delta)$, forming a subset touching at least $\delta n$ edges -- a contradiction.)

The above two propositions can -- and will -- be used to argue that a sparse random graph contains typically a linearly sized locally sparse subgraph of bounded maximum degree.

We can now formulate the main result of the paper.

\begin{thm}\label{thm1}
Let $c_1>c_2>1$, $0<\alpha<1$, $\Delta>0$. Let $G=(V,E)$ be a graph on  $n$ vertices, satisfying:
\begin{enumerate}
\item $\frac{|E|}{|V|}\ge c_1$\,;
\item every vertex subset $W\subset V$ of size $|W|\le \alpha n$ spans less than $c_2|W|$ edges;
\item $\Delta(G)\le \Delta$.
\end{enumerate}
Then $G$ contains an induced subgraph $G^*=(V^*,E^*)$ on at least $\alpha n$ vertices which is a $\gamma$-expander, for $\gamma=\frac{c_1-c_2}{\Delta\cdot\left\lceil\log_2\frac{1}{\alpha}\right\rceil}$.
\end{thm}

Putting it informally, every locally sparse graph $G$ of bounded maximum degree contains a linearly sized expander $G^*$. In our terminology, the first two conditions above say precisely that $G$ is a $(c_1,c_2,\alpha)$-graph. They are spelled out in full in the statement above so as to make it self-contained.

With some sacrifice in constants involved (and perhaps in transparency of the proof) we can make our argument algorithmic. The proof relies on the Cheeger inequality, providing a fairly standard nowadays connection between graph eigenvalues and expansion. The obtained result is summarized in the following theorem.
\begin{thm}\label{thm2}
Let $c_1>c_2>1$, $0<\alpha<1$, $\Delta>0$. There exist a constant $\gamma=\gamma(c_1,c_2,\alpha,\Delta)>0$ and an algorithm, that, given an $n$-vertex graph $G=(V,E)$ with $\Delta(G)\le \Delta$ and $|E|/|V|\ge c_1$, finds in time polynomial in $n$ a subset $W\subset V$ of size $|W|\le \alpha n$, spanning at least $c_2|W|$ edges in $G$, or an induced $\gamma$-expander $G^*=(V^*,E^*)\subseteq G$ on at least $\alpha n$ vertices.
\end{thm}

Let us now discuss the result and compare it with prior results about finding expanders in given graphs. Observe first that not every bounded degree graph $G$ of density $|E(G)|/|V(G)|=c_1>1$ contains a linearly sized expander. For example, if $G$ is taken to be the $\sqrt{n}\times\sqrt{n}$ grid, then $G$ is a planar graph of density close to 2. Since every planar graph $G^*$ on $m$ vertices has a separator of size $O(\sqrt{m})$ by the famous Lipton--Tarjan theorem \cite{LT79}, and an $m$-vertex $\gamma$-expander has all its separators of linear size (see Section \ref{sec-minors} for discussion on separators), the graph $G$ does not contain a linearly sized expander. We can conclude that the local sparsity condition is necessary. The assumption about bounded maximum degree of $G$ will only be used in the proof to convert edge expansion into vertex expansion -- we will argue that in the final graph $G^*=(V^*,E^*)$ for every subset $W\subset V^*$, $|W|\le |V^*|/2$, $G^*$  has many edges crossing between $W$ and its complement. Then the assumption $\Delta(G)\le \Delta$ will imply large vertex boundary of such $W$ as well.

As we indicated above, Shapira and Sudakov \cite{SS15} and later Montgomery \cite{Mon15} argued that {\em every} graph $G$ contains a weak expander $G^*$ of nearly the same density. Their notion of expansion is different -- the expansion required is gradual in terms of subset sizes; for an $m$-vertex graph to be a weak expander, vertex sets of size $m^c$, $0<c<1$, should expand by a constant factor, whereas linearly sized vertex sets are required to expand only by about $1/\log m$ factor. (We are rather informal here in our descriptions, see the actual papers \cite{SS15,Mon15} for accurate definitions.) Neither of these results guarantees a  (weakly) expanding subgraph on linearly many vertices. One should probably add here that both papers \cite{SS15,Mon15} are mainly concerned not with expanders, but rather with finding complete minors of relatively small order; the expansion statements serve as a technical tool and are pretty much tailored to that particular target. Moreover, Moshkovitz and Shapira \cite{MS17} argue that the requirement $|N_{G^*}(W)|=\Omega\left({|W|}/{\log m}\right)$ cannot be strengthened much; their example of a graph $G$ showing it is a bounded degree graph of density $c_1>1$, but again is not locally sparse. Much earlier, Koml\'os and Szemer\'edi, in their work on topological cliques in graphs \cite{KS94,KS96}, presented a fairly general scheme for arguing about the existence of weak expanders in any given graph; their scheme does not provide -- naturally -- for finding linearly sized expanders, or expanders with constant expansion of subsets.

The rest of the paper is structured as follows. In the next section we prove Theorems \ref{thm1} and \ref{thm2}.
In Section \ref{sec-random} we discuss random graphs and in particular prove Propositions \ref{prop1} and \ref{prop2}. In Section \ref{sec-minors} we present statements about embedding large minors in expanding graphs, locally sparse graphs and random graphs. In Section \ref{sec-games} we discuss results about positional games, namely, about the biased minor creation games in their Maker--Breaker, Avoider--Enforcer and Client--Waiter versions.

\medskip

\noindent{\bf Notation.} Our notation is mostly standard. As stated before, for a graph $G=(V,E)$ and a vertex subset $W\subset V$ we denote by $N_G(W)$ the external neighborhood of $W$ in $G$. If $U,W$ are disjoint subsets of $G$, then $e_G(U,W)$ stands for the number of edges of $G$ between $U$ and $W$. We suppress the rounding notation occasionally to simplify the presentation.

\section{Proof of the main result}\label{sec-main}
We first describe briefly (and rather informally) the main idea of the proof of Theorem \ref{thm1}. The argument proceeds iteratively, initiating with the given graph $G$. At a current iteration $i$, if the current graph $G_i$ contains a subset $W$ of size $|W|\le |V(G_i)|/2$ whose density is close enough to that of $G_i$, we update $G_{i+1}:=G[W]$. Otherwise we argue that $G_i$ (more accurately, a subgraph of it -- see the forthcoming proof for details) is a good edge expander, which can be translated immediately to good vertex expansion due to the assumption of the bounded maximum degree of $G$. This strategy is somewhat reminiscent of the approach taken in \cite{SS15, Mon15}.

Let us get to the actual proof. Set
$$
\delta= \frac{c_1-c_2}{\left\lceil\log_2\frac{1}{\alpha}\right\rceil}\,;\quad\quad \beta_i=c_1-i\delta\,,\quad i\ge 0\,.
$$
We initialize with $i=0$, $G_0=G$.

Suppose now we are at iteration $i\ge 0$, and the current graph $G_i$ satisfies $|V(G_i)|\le n/2^i$, $|E(G_i)|/|V(G_i)|\ge \beta_i$. (This is obviously true for $i=0$.) Let $H_i=(U_i,F_i)$ be a minimal by inclusion non-empty induced subgraph of $G_i$  for which $|F_i|/|U_i|\ge \beta_i$. (Such a subgraph exists, as $G_i$ itself meets the requirements.) Then every subset $W\subset U_i$ touches at least $\beta_i|W|$ edges of $H_i$. Otherwise, deleting $W$ is easily seen to produce a smaller induced subgraph $H_i'$, still meeting the requirement stated in the definition of $H_i$ -- a contradiction. If
\begin{equation}\label{update_rule}
\mbox{there exists $W\subset U_i$, $\alpha n\le |W|\le \frac{|U_i|}{2}$\,, s.t. $W$ spans at least $\beta_{i+1}|W|$ edges in $H_i$}\,,
\end{equation}
then we update $G_{i+1}:=G[W]$, $i:=i+1$, else we halt the process.

Observe that if we indeed proceed to the next iteration as described above, then the set $W$ used to create the next graph is at most half of $V_i$ in size, and has density at least $\beta_{i+1}$ by (\ref{update_rule}), so the new graph $G_{i+1}=G[W]$ satisfies the required inductive assumptions.

If the above defined iterative process gets to $i=\left\lceil\log_2\frac{1}{\alpha}\right\rceil$, then the graph $G_i$ satisfies: $|V(G_i)|\le n/2^i\le \alpha n$ and $|E(G_i)|/|V(G_i)|\ge \beta_i=c_2$ -- a contradiction to our assumption on $G$, which is postulated to be locally sparse. Hence the process stops with $i\le \left\lceil\log_2\frac{1}{\alpha}\right\rceil-1$, with the reason being that the iteration condition (\ref{update_rule}) is not met in the current/last graph $H_i=(U_i,F_i)$. The density of $H_i$ is at least $\beta_i\ge c_2+\delta$, implying in particular that $H_i$ has at least $\alpha n$ vertices.

Now, if $W\subset U_i$ has size $\alpha n\le |W|\le |U_i|/2$, then upon looking at (\ref{update_rule}) we conclude that $W$ spans at most $\beta_{i+1}|W|$ edges; yet, due to the choice of $H_i$, the same $W$ touches at least $\beta_i|W|$ edges in $H_i$. It follows that $H_i$ has at least $(\beta_i-\beta_{i+1})|W|=\delta|W|$ edges crossing between $W$ and its complement $U_i\setminus W$. For the complementary case $|W|\le \alpha n$, by the local density assumption we get that $W$ spans at most $c_2|W|$ edges, yet touches at least $\beta_i|W|\ge (c_2+\delta)|W|$ edges in $H_i$. This again implies that $H_i$ has at least $\delta|W|$ edges crossing between $W$ and its complement $U_i\setminus W$. So the edge boundary of $W$ in $H_i$ is at least $\delta |W|$, and recalling the assumption $\Delta(G)\le \Delta$, we derive that $|N_{H_i}(W)|\ge \frac{\delta |W|}{\Delta}=\frac{c_1-c_2}{\Delta\cdot\left\lceil\log_2\frac{1}{\alpha}\right\rceil}|W|$. Hence $H_i$ meets the required specifications for a sought $\gamma$-expander, and we can take $G^*=H_i$ and complete the proof.\hfill$\Box$

\bigskip

We now discuss the algorithmic aspect of the problem and the proof of Theorem \ref{thm2}. (Thanks to Noga Alon and Avi Wigderson for bringing up the algorithmic issue, and for their input.) First, the required background briefly. (See \cite{Chung-book} for an extensive discussion of the subject and all missing definitions.) For a graph $G=(V,E)$ on $n$ vertices and a subset $W\subseteq V$ set $vol_G(W)=\sum_{v\in W}deg_G(v)$. Let $0=\lambda_0\le \lambda_1\le\ldots\le \lambda_{n-1}$ be the eigenvalues of the normalized Laplacian of $G$, ordered in the non-decreasing order, and set $\lambda(G)=\lambda_1$. Obviously $\lambda(G)$ can be computed in time polynomial in $n$.

Let now
$$
h(G)= \min_{\emptyset\ne W\subsetneq V}\, \frac{e_G(W,V\setminus V)}{\min(vol_G(W),vol_G(V\setminus W))}\,,
$$
the quantity $h(G)$ is sometimes called the {\em Cheeger constant} of $G$. Having $h(G)$ large means the graph is a good edge expander, and -- assuming its maximum degree is bounded (this is the assumption we adopt throughout this discussion) -- is a good vertex expander.  The famous {\em Cheeger inequality for graphs} \cite{Dod84, AM85, Alo86} states that
\begin{equation}\label{Cheeger}
\frac{h^2(G)}{2}\le \lambda(G)\le 2h(G)\,.
\end{equation}
(This is one of the most convincing manifestations of the famous and very fruitful connection between expansion of a graph and its eigenvalues.)

The proof of the first inequality in (\ref{Cheeger}) is constructive in the sense that it finds in time polynomial in $n$ a subset $W\subset V$, $vol_G(W)\le vol_G(V)/2$, satisfying $e_G(W,V\setminus W)\le \sqrt{2\lambda(G)}\cdot vol_G(W)$. (See \cite{Alo86}, or \cite[Ch. 2]{Chung-book}, or \cite{Chu10}, or \cite[Sect. 4.5]{HLW06}.)

We can now describe (somewhat informally) an algorithm for finding a large expander in an input graph, or discovering a small and dense subset of vertices in it. Let $c_1>c_2>1$, $0<\alpha<1$, $\Delta>0$ be fixed parameters. The algorithm gets an $n$-vertex graph $G=(V,E)$, satisfying $|E|/|V|\ge c_1$, $\Delta(G)\le \Delta$, as an input. Similarly to the proof of Theorem \ref{thm1}, the algorithm proceeds in iterations, starting from $V_0=V$, $i=0$. Suppose we are at iteration $i$ with $G_i=G[V_i]=(V_i,E_i)$ as a current graph of density $d_i=|E_i|/|V_i|\ge c_2$. In case $|V_i|\le \alpha n$ we are done -- a small and dense subset has been found. Assume otherwise. If $G_i$ has an isolated vertex $v$, we update $V_{i+1}:=V_i\setminus\{v\}$, $i:=i+1$, and proceed to the next iteration. Otherwise,
compute the eigenvalue $\lambda(G_i)$. If $\lambda(G_i)$ is large, then by the second part of the Cheeger inequality (\ref{Cheeger}) we obtain that the current graph $G_i$ is a good edge expander, and is thus a good vertex expander due to the assumption $\Delta(G)\le \Delta$. If $\lambda(G_i)$ is small, then the algorithm finds in time polynomial in $n$ a subset $W_i\subset V_i$ with $vol_{G_i}(W_i)\le vol_{G_i}(V_i)/2$ and small edge boundary, say, $e_{G_i}(W_i,V_i\setminus W_i)\le \delta |W_i|$. Observe that $\Delta\cdot |V_i\setminus W_i|\ge vol_{G_i}(V_i\setminus W_i)\ge \frac{vol_{G_i}(V_i)}{2}\ge \frac{|V_i|}{2}$, implying $|W_i|\le \left(1-\frac{1}{2\Delta}\right)|V_i|$.
If the set $W_i$ touches at most $d_i|W_i|$ edges in $G_i$, the algorithm deletes $W_i$ and updates $V_{i+1}:=V_i\setminus W_i$, $i:=i+1$. If $W_i$ touches at least $d_i|W_i|$ edges in $G_i$, we obtain $e_{G_i}(W_i)\ge (d_i-\delta)|W_i|$. We then update $V_{i+1}:=W_i$, $i:=i+1$. Finally, we proceed to the next iteration.
Now, to put all things together we can set $k=\frac{\log\alpha}{\log\left(1-\frac{1}{2\Delta}\right)}$, $\delta=\frac{c_1-c_2}{2k}=\frac{(c_1-c_2)\log\frac{2\Delta}{2\Delta-1}}{2\log\frac{1}{\alpha}}$. This completes (the sketch of) the proof of Theorem \ref{thm2}. 

(Since the above algorithmic argument loses to the existential proof of Theorem \ref{thm1} both in terms of the bound delivered and of transparency, we allowed ourselves to be somewhat informal in the algorithmic description above).

\section{Random graphs}\label{sec-random}
We first prove Propositions \ref{prop1} and \ref{prop2} here. Both proofs are pretty straightforward; in fact, Proposition \ref{prop1} is stated and proven in \cite{Kri16}, we reproduce its proof here for the sake of completeness.

\noindent{\bf Proof of Proposition \ref{prop1}.} The probability in $G(n,c_1/n)$ that there exists a vertex subset violating the required property is at most
\begin{eqnarray*}
&&\sum_{i\le \alpha n}\binom{n}{i}\binom{\binom{i}{2}}{c_2i}\cdot p^{c_2i}\le
 \sum_{i\le \alpha n}\left(\frac{en}{i}\right)^i\cdot \left(\frac{eip}{2c_2}\right)^{c_2i}
=\sum_{i\le\alpha n}\left[\frac{en}{i}\cdot\left(\frac{ec_1i}{2c_2n}\right)^{c_2}\right]^i\\
&=&\sum_{i\le\alpha n}\left[\frac{e^{c_2+1}c_1^{c_2}}{(2c_2)^{c_2}}\cdot\left(\frac{i}{n}\right)^{c_2-1}\right]^i\,.
\end{eqnarray*}

Denote the $i$-th summand of the last sum by $a_i$. Then, if $i\le n^{1/2}$ we get: $a_i\le \left(O(1)n^{-\frac{c_2-1}{2}}\right)^i$, implying $\sum_{i\le n^{1/2}}a_i=o(1)$. For $n^{1/2}\le i\le \alpha n$, we have, recalling the expression for $\alpha$:
$$
a_i\le \left[\frac{e^{c_2+1}c_1^{c_2}}{(2c_2)^{c_2}}\cdot\left(\frac{c_2}{5c_1}\right)^{c_2}\right]^i\le \left(e\cdot\left(\frac{e}{10}\right)^{c_2}\right)^i=o(n^{-1/2})\,.
$$
It follows that $\sum_{i\le\alpha n}a_i=o(1)$, and the desired property of the random graphs holds {\bf whp}. \hfill$\Box$

\bigskip

\noindent{\bf Proof of Proposition \ref{prop2}.} The probability in $G(n,C/n)$ that there exists a vertex subset violating the required property is at most
$$
\binom{n}{\frac{\delta n}{\ln\frac{1}{\delta}}}\binom{\frac{\delta n^2}{\ln\frac{1}{\delta}}}{\delta n} \left(\frac{C}{n}\right)^{\delta n}
= \left[ \left(\frac{e\ln\frac{1}{\delta}}{\delta}\right)^{\frac{1}{\ln\frac{1}{\delta}}}\cdot \frac{eC}{\ln\frac{1}{\delta}}\right]^{\delta n}\,.
$$
Taking $\delta=\delta(C)$ to be small enough guarantees that the expression above vanishes with growing $n$. \hfill$\Box$

\bigskip

We now use our main result to argue that a supercritical random graph $G(n,c/n)$, $c>1$, contains \whp an induced expander of linear size.

\begin{corol}\label{cor1}
For every $\epsilon>0$ there exists $\gamma>0$ such that a random graph $G\sim G\left(n,\frac{1+\epsilon}{n}\right)$ contains \whp an induced bounded degree  $\gamma$-expander on at least $\gamma n$ vertices.
\end{corol}

\noindent{\bf Proof.} Due to the standard monotonicity arguments we can assume that $\epsilon$ is small enough where necessary.

We will utilize several (very standard) facts about supercritical random graphs. It is known (see, e.g.,
\cite[Ch. 5]{JLR}) that \whp $G\sim G\left(n,\frac{1+\epsilon}{n}\right)$ contains a connected component $C_1=(V_1,E_1)$ (the so called giant component) satisfying:
\begin{eqnarray*}
|V_1|&=&2\epsilon (1+o_{\epsilon}(1))n\,,\\
\frac{|E_1|}{|V_1|} &=& 1+(1+o_{\epsilon}(1))\frac{\epsilon^2}{3}\,.
\end{eqnarray*}

Also, by Proposition \ref{prop2} \whp every $\frac{\epsilon^3}{2\ln\frac{1}{\epsilon}}n$ vertices of $C_1$ touch at most $\frac{\epsilon^3}{3}n$ edges. Deleting $\frac{\epsilon^3}{2\ln\frac{1}{\epsilon}}n$ vertices of highest degrees from  $C_1$, one gets a graph $G_0=(V_0,E_0)$ of maximum degree $\Delta(G_0)\le 4\ln\frac{1}{\epsilon}$. In addition,
\begin{eqnarray*}
|V_0|&\ge &\left(2\epsilon (1+o_{\epsilon}(1))-\frac{\epsilon^3}{2\ln\frac{1}{\epsilon}}\right)n
=2\epsilon (1+o_{\epsilon}(1))n\,,\\
|E_0|&\ge& |E_1|-\frac{\epsilon^3}{3}n
\ge |V_1|\left(1+\frac{\epsilon^2}{3}+o(\epsilon^2)\right)-\frac{\epsilon^3}{3}n\\
 &\ge& |V_0|\left(1+\frac{\epsilon^2}{3}+o(\epsilon^2)\right)-\frac{\epsilon^3}{3}n\ge |V_0|\left(1+\frac{\epsilon^2}{7}\right).
\end{eqnarray*}
Finally, applying Proposition \ref{prop1} with $c_1=1+\epsilon$, $c_2=1+\frac{\epsilon^2}{10}$ we get that \whp every $k\le \alpha n$ vertices of $G_0$ (with $\alpha=\alpha(\epsilon)$ from Proposition \ref{prop1}) span fewer than $(1+\epsilon^2/10)k$ edges. The conditions are set to call Theorem \ref{thm1} and to apply it to $G_0$; we conclude that, given the above likely events, $G_0$ contains a linearly sized $\gamma$-expander.  \hfill$\Box$

\bigskip

We remark here that in order to get the above stated qualitative result, one does not really need to apply the heavy machinery of random graphs -- it is enough actually to argue from the ``first principles". Indeed, given the likely existence of the giant component $C_1$ in the supercritical regime, one can argue easily (for example, through sprinkling) that its density is typically above 1. In quantitative terms, the above argument delivers very weak (but constant) expansion; much more accurate results can be obtained by invoking a powerful statement of Ding, Lubetzky and Peres \cite{DLP14}, describing in great detail a likely structure of the giant component in the  supercritical regime. Working out the details carefully, based on \cite{DLP14}, should probably deliver the likely existence of a $\Theta(\epsilon)$-expander on $\Theta(\epsilon)n$ vertices. We chose not to perform this (pretty technical) analysis here, relying instead on our main theorem to get a qualitative result quickly.

\section{Separators and embedding minors}\label{sec-minors}
In this section we discuss the relation between graph expansion and graph separators, and also argue that expanders -- and thus locally sparse graphs -- contain large minors.

Given a graph $G=(V,E)$ on $n$ vertices, a vertex set $S\subset V$ is called a {\em separator} if there is a partition $V=A\cup B\cup S$ of the vertex set of $G$ such that $G$ has no edges between $A$ and $B$, and $|A|,|B|\le 2n/3$. Separators serve to measure quantitatively the connectivity of large vertex sets in graphs; the fact that all separators in $G$ are large indicates that it is costly to break $G$ into large pieces not connected by any edge.

It is easy to argue that expanders do not have small separators. Indeed, let $G=(V,E)$ be a $\gamma$-expander on $n$ vertices, and let $S$ be a separator in $G$ of size $|S|=s$, separating $A$ and $B$, with $|A|=a$, $|B|=b$; we assume $a\le b\le 2n/3$. Then $a+s\ge n/3$. Clearly, $N_G(A)\subseteq S$. Since $a\le n/2$, by the definition of a $\gamma$-expander we get $s-\gamma a\ge 0$. Multiplying this inequality by $1/\gamma$ and summing with $a+s\ge n/3$, we obtain: $s(1+1/\gamma)\ge n/3$, or $s\ge \frac{\gamma n}{3(\gamma+1)}$. We have proven:

\begin{propos}\label{prop3}
Let $G$ be a $\gamma$-expander on $n$ vertices, and let $S$ be a separator in $G$. Then  $|S|\ge \frac{\gamma n}{3(\gamma+1)}$.
\end{propos}

We now discuss embedding minors in expanders and in locally sparse graphs. Let $G=(V,E)$, $H=(U,F)$ be graphs with $U=\{u_1,\ldots, u_t\}$. We say that $G$ contains $H$ as a {\em minor} if there is a collection $(V_1,\ldots,V_t)$ of pairwise disjoint vertex subsets in $V$ such that each $V_i$ spans a connected subgraph in $G$, and whenever $(u_i,u_j)\in F$, the graph $G$ has an edge between $V_i$ and $V_j$. (Then contracting each $U_i$ to a single vertex produces a copy of $H$.) Minors are one of the most important concepts in graph theory, and finding sufficient conditions for embedding minors is a very frequently considered research direction. Observe trivially that if $G$ contains a minor of $H$ then $|V(G)|\ge |V(H)|$ and $|E(G)|\ge |E(H)|$; these trivial bounds provide an obvious but meaningful benchmark for minor embedding statements.

Kleinberg and Rubinfeld proved in \cite{KR96} that a $\gamma$-expander on $n$ vertices  of maximum degree $\Delta$ contains all graphs with $O(n/\log^{\kappa}n)$ vertices and edges as minors, for $\kappa=\kappa(\gamma,\Delta)>0$. Phrasing it differently, a bounded degree $\gamma$-expander is {\em minor universal} for all graphs with $O(n/\log^{\kappa}n)$ vertices and edges. This is optimal up to polylogarithmic factors due to the above stated trivial bound, as there exist $n$-vertex expanders with $\Theta(n)$ edges. (Formally, the number of vertices $n$ is another bottleneck here.)  From Theorem \ref{thm1} we obtain the following corollary.

\begin{corol}\label{minors}
For every $c_1>c_2>1$, $0<\alpha<1$, $\Delta>0$ there exists $\kappa>0$ such that every $(c_1,c_2,\alpha)$-graph $G$ on $n$ vertices of maximum degree at most $\Delta$ contains all graphs with at most $n/\log^{\kappa}n$ vertices and edges as minors.
\end{corol}

So in particular, recalling Corollary \ref{cor1}, we conclude that a supercritical random graph $G\sim G\left(n,\frac{1+\epsilon}{n}\right)$ is \whp minor-universal for the family of graphs with at most $n/\log^{\kappa}n$ vertices and edges, for some $\kappa=\kappa(\epsilon)>0$. Proving that $\kappa$ in the above statement can be taken to be independent of $\epsilon$ is an open problem; establishing the best possible value of $\kappa$ would be quite nice.

We now switch to embedding complete minors in locally sparse graphs. Kawarabayshi and Reed, completing a long and illustrious line of research in this topic (see \cite{LT79, AST90, PRS94} for some milestones) proved in \cite{KR10} that a graph $G$ on $n$ vertices has a minor of the complete graph $K_h$ or a separator of order $O(h\sqrt{n})$. Since a $\gamma$-expander $G$ on $n$ vertices has all separators linear in size by Proposition \ref{prop3}, we conclude from Theorem \ref{thm1}:

 \begin{corol}\label{complete_minors}
For every $c_1>c_2>1$, $0<\alpha<1$, $\Delta>0$ there exists $c>0$ such that every $(c_1,c_2,\alpha)$-graph $G$ on $n$ vertices of maximum degree at most $\Delta$ contains a minor of $K_{c\sqrt{n}}$.
\end{corol}

As due to Corollary \ref{cor1} a supercritical random graph $G\sim G\left(n,\frac{1+\epsilon}{n}\right)$ contains \whp a bounded degree $\gamma$-expander on $\Theta(n)$ vertices, we obtain that such a random graph typically contains a minor of the complete graph on $\Theta(\sqrt{n})$ vertices. This recovers the result of Fountoulakis, K\"uhn and Osthus \cite{FKO08}, obtained through direct (and quite involved technically) means.

\section{Positional games}\label{sec-games}
We will now discuss applications of our main result to positional games (see \cite{HKSS-book} for a systematic introduction to this fascinating combinatorial discipline). More specifically, we will address  minor creation games. The game types we will cover are Maker--Breaker games, Avoider--Enforcer games and Client--Waiter games. All games to be considered are played on the edge set of the complete graph $K_n$ on $n$ vertices; we assume the parameter $n$ to be sufficiently large where necessary.

Our arguments for all three game types make use of the following two families of subgraphs of $K_n$, parameterized by $\epsilon,\delta$: let ${\cal F}_1={\cal F}_1(\epsilon,\delta)$ be the family of all subgraphs of $K_n$ with $\frac{\epsilon n}{4}$ edges and covering number at most $\delta n$; let ${\cal F}_2={\cal F}_2(\epsilon,\delta)$ be the family of all subgraphs $H=(U,F)$ on at most $\delta n$ vertices and of density $|F|/|U|= 1+\epsilon/8$.

\noindent{\bf Maker--Breaker games.} In a {\em Maker--Breaker} game two players, called Maker and Breaker, claim alternately free edges of the complete graph $K_n$, with Maker moving first. Maker claims one edge at a time, while Breaker claims $b\ge 1$ edges (or all remaining fewer than $b$ edges if this is the last round of the game). The integer parameter $b$ is the so-called {\em game bias}. Maker wins the game if the graph of his edges in the end possesses a given graph theoretic property, Breaker wins otherwise, with draw being impossible. In the minor creation game Maker's goal is to create a minor of the complete graph $K_t$ for $t=t(n)$ as large as possible. This game has been considered by Hefetz, Krivelevich, Stojakovi\'c and Szab\'o in \cite{HKSS08}. They established a kind of a sharp phase transition at $b=n/2$. For $b\ge n/2$, as follows from a general result by Bednarska and Pikhurko \cite{BP05}, Breaker has a strategy to force Maker's graph being acyclic (and thus not containing a $K_3$-minor) by the end of the game; for $b=(1-\epsilon)n/2$ \cite{HKSS08} showed that Maker has a strategy to create a complete minor of order $c\sqrt{n/\log n}$ for $c=c(\epsilon)>0$. Here we improve the latter result to the optimal order of magnitude by proving:

\begin{thm}\label{thm-MB}
For every $\epsilon>0$ there exists $c>0$ such that for all sufficiently large $n$, when playing a $b$-biased Maker--Breaker game on $E(K_n)$ with $b\le (1-\epsilon)\frac{n}{2}$, Maker has a strategy to create a minor of the complete graph $K_{c\sqrt{n}}$.
\end{thm}

\noindent{\bf Proof.} Due to the bias monotonicity we can assume that $\epsilon$ is small enough where necessary. Maker's strategy is very simple: during the first $(1+\frac{\epsilon}{2})n$ rounds he plays {\em randomly}, i.e., chooses a uniformly random edge to claim out of the set of all available edges at that moment. Observe that in each of these rounds, the number of available edges is at least
$$
\binom{n}{2}-(b+1)\left(1+\frac{\epsilon}{2}\right)n\ge \frac{\epsilon}{3}n^2\,,
$$
meaning that the probability an edge $e$ is chosen by Maker is at most $3/(\epsilon n^2)$, regardless of the history of the game. Therefore, for a subgraph $H$ of $K_n$ the probability that the graph $M_0$ of Maker's edges after the first $(1+\frac{\epsilon}{2})n$ rounds contains $H$ is at most:
$$
\left(\left(1+\frac{\epsilon}{2}\right)n\right)^{|E(H)|}
\left(\frac{1}{\frac{\epsilon}{3}n^2}\right)^{|E(H)|}
= \left(\frac{3\left(1+\frac{\epsilon}{2}\right)}{\epsilon n}\right)^{|E(H)|}
\le\left(\frac{6}{\epsilon n}\right)^{|E(H)|}
$$
(for each edge $e\in E(H)$ decide in which of the first $\left(1+\frac{\epsilon}{2}\right)n$ rounds $e$ is to be claimed by Maker, and require that all  edges are indeed claimed in the chosen rounds). It follows that the probability that $M_0$ contains any graph $H$ from ${\cal F}_1$ can be bounded from above by:
$$
\binom{n}{\delta n}\binom{\delta n^2}{\frac{\epsilon n}{4}}\cdot\left(\frac{6}{\epsilon n}\right)^{\frac{\epsilon n}{4}}
= \left[\left(\frac{e}{\delta}\right)^{\delta}\cdot \left(\frac{24e\delta}{\epsilon^2}\right)^{\frac{\epsilon}{4}}\right]^n
= o(1)\,,
$$
for $\delta\le\delta_1$ with $\delta_1=\delta_1(\epsilon)$ small enough. The probability that $M_0$ contains any graph $H$ from ${\cal F}_2$ can be bounded from above by:
$$
\sum_{k\le \delta n}\binom{n}{k}\binom{\binom{k}{2}}{\left(1+\frac{\epsilon}{8}\right)k}\left(\frac{6}{\epsilon n}\right)^{\left(1+\frac{\epsilon}{8}\right)k}
\le \sum_{k\le \delta n}\left[\frac{en}{k}\cdot\left(\frac{3ek}{\epsilon n}\right)^{1+\frac{\epsilon}{8}}\right]^k
\le \sum_{k\le \delta n}\left[ \frac{10e^2}{\epsilon^2}\cdot \left(\frac{k}{n}\right)^{\frac{\epsilon}{8}}\right]^k
= o(1)\,,
$$
for $\delta\le\delta_2$ with $\delta_2=\delta_2(\epsilon)$ small enough. Take $\delta=\min\{\delta_1,\delta_2\}$, then with positive probability Maker, playing against any strategy of Breaker, can create a graph $M_0$ on $n$ vertices with the following properties:
\begin{itemize}
\item[{\bf (P1)}] has at least $\left(1+\frac{\epsilon}{2}\right)n$ edges;
\item[{\bf (P2)}] every $k\le\delta n$ vertices span at most $\left(1+\frac{\epsilon}{8}\right)k$ edges;
\item[{\bf (P3)}] every $\delta n$ vertices touch at most $\frac{\epsilon n}{4}$ edges.
\end{itemize}
Since the game analyzed is a perfect information game with no chance moves, it follows that in fact Maker has a (deterministic) strategy to create a graph $M_0$ with the above stated properties. Take such $M_0$ and delete $\delta n$ vertices of highest degrees. The obtained graph $M_1$ has $(1-\delta)n$ vertices, at least $\left(1+\frac{\epsilon}{4}\right)n$ edges, maximum degree $\Delta(M_1)\le \frac{\epsilon}{2\delta}$, and every $k\le \delta n$ vertices span at most $\left(1+\frac{\epsilon}{8}\right)k$ edges. Applying Corollary \ref{complete_minors} shows that such $M_1$, being a part of Maker's graph by the end of the game, contains a complete minor $K_{\Theta(\sqrt{n})}$. \hfill$\Box$

\bigskip

\noindent{\bf Avoider--Enforcer games.} In a biased {\em Avoider--Enforcer} game two players, called Avoider and Enforcer, claim alternately free edges of the complete graph $K_n$; for simplicity we assume Enforcer moves first (this will not change much for the games to be considered). Avoider claims one edge at a time, while Enforcer claims exactly $b\ge 1$ edges, or all remaining fewer than $b$ edges if this is the last round of the game. (The rules we describe here are the so called strict rules, there is also the monotone version of the rules, see, e.g., \cite{HKSS10} for discussion.) Avoider--Enforcer games are a mis\'ere version of Maker--Breaker games -- Avoider wins the game if the graph of his edges in the end {\em does not} possess a given graph theoretic property, Enforcer wins otherwise, with draw being impossible. In the minor creation game Enforcer's goal is to force a minor of the complete graph $K_t$ in Avoider's final graph, for $t=t(n)$ as large as possible. This game has too been considered by Hefetz, Krivelevich, Stojakovi\'c and Szab\'o in \cite{HKSS08}. They showed that for $b\le n/19$ Enforcer can force a minor of $K_{c\sqrt{n/\log n}}$ in Avoider's graph, and that for $b<(1-\epsilon)n/2$, $\epsilon>0$ a constant, Enforcer can force a complete minor of order $n^{\delta}$ for $\delta=\delta(\epsilon)>0$. The proof idea of the latter result in \cite{HKSS08} was quite different from the present argument: they argued that Enforcer can force Avoider to create a graph $A$ on $n$ vertices with about $(1+\epsilon)n$ edges and with $o(n)$  cycles of length $O(\log n)$; deleting one edge from each short cycle produces a graph $A_1$ on $n$ vertices with at least $(1+\epsilon/2)n$ edges and of girth $\Omega(\log n)$. Then, Hefetz et al. showed in \cite{HKSS08}, invoking a result of K\"uhn and Osthus from \cite{KO03}, that such a graph contains  a polynomially large complete minor.

We provide an improvement of the above results to the optimal order of magnitude by proving:

\begin{thm}\label{thm-AE}
For every small enough $\epsilon>0$ there exists $c>0$ such that for all sufficiently large $n$, when playing a $b$-biased Avoider--Enforcer game on $E(K_n)$ with $b\le (1-\epsilon)\frac{n}{2}$, Enforcer has a strategy to force a minor of the complete graph $K_{c\sqrt{n}}$ in Avoider's graph.
\end{thm}

\noindent{\bf Proof.} Enforcer's goal is to put his edge (or to break into) every graph in the family ${\cal F}_1\cup{\cal F}_2$, for some small $\delta$ from the definition of ${\cal F}_1,{\cal F}_2$ to be set later. He thus disguises himself as Breaker and uses the following criterion for Breaker's win due to Beck \cite{Bec82}. (We present it here in the form adapted for our setting.)

\begin{lemma}\cite{Bec82}\label{genES}
Let $n,b$ be positive integers, and let ${\cal F}$ be a family of subgraphs of $K_n$. If
$$
\sum_{H\in {\cal F}} (1+b)^{-|E(H)|}<1\,,
$$
then in the $b$-biased Maker--Breaker game on $E(K_n)$, Breaker, as the first player to move, has a strategy to put his edge into every $H\in {\cal F}$.
\end{lemma}

Now we need to crunch some numbers. For the family ${\cal F}_1$ we have:
$$
\sum_{H\in{\cal F}_1} (1+b)^{-|E(H)|}\le \binom{n}{\delta n}\binom{\delta n^2}{\frac{\epsilon n}{4}}\left(\frac{1}{1+(1-\epsilon)\frac{n}{2}}\right)^{\frac{\epsilon n}{4}}
\le \left[\left(\frac{e}{\delta}\right)^{\delta}\cdot \left(\frac{30\delta}{\epsilon}\right)^{\frac{\epsilon}{4}}\right]^n=o(1)\,,
$$
for $\delta\le\delta_1$. For the family ${\cal F}_2$ the calculation gives:
$$
\sum_{H\in{\cal F}_1} (1+b)^{-|E(H)|}\le \sum_{k\le \delta n}\binom{n}{k}\binom{\binom{k}{2}}{\left(1+\frac{\epsilon}{8}\right)k}
\left(\frac{1}{1+(1-\epsilon)\frac{n}{2}}\right)^{\left(1+\frac{\epsilon n}{8}\right)k}=o(1)\,,
$$
for $\delta\le\delta_2$. Taking $\delta=\min\{\delta_1,\delta_2\}$ and denoting ${\cal F}={\cal F}_1\cup{\cal F}_2$, we see that the condition of Lemma \ref{genES} holds, and hence Enforcer can put his edge into every graph $H\in {\cal F}$. It follows that Avoider's graph in the end of the game satisfies conditions {\bf (P1)}--{\bf (P3)} above, and arguing in essentially the same way as in the proof of Theorem \ref{thm-MB}, we conclude that Enforcer has a strategy to force a complete minor $K_t$ with $t=\Theta(\sqrt{n})$ in Avoider's graph. \hfill$\Box$

\bigskip

\noindent{\bf Client--Waiter games.} In a {\em Client--Waiter} game with bias $b$ there are two players, called Client and Waiter. The game proceeds in rounds, where in each round of the game Waiter offers to Client between one and $b+1$ edges of $K_n$, previously not offered by him. Client claims an edge of his choice among the edges offered, and the remaining edges are assigned to Waiter. The game runs till all edges of $K_n$ have been offered. Client wins the game if his final graph, composed of all the edges claimed by him during the game, possesses a target graph theoretic property, Waiter wins otherwise, with draw being impossible. Client--Waiter games, also called Chooser--Picker games by Beck (see, e.g., \cite{Bec02}), have been quite popular in recent research. In the Client--Waiter minor game Client aims to create as large a complete minor as possible. This game has been considered by Hefetz, Tan and the author in \cite{HKT16}. The authors of \cite{HKT16} observed (based on a result from \cite{BHKL16}) that for $b\ge n/2-1$ Waiter has a strategy to keep Client's graph acyclic (and thus $K_3$-minor-free) throughout the game. On the other hand, they proved that for any fixed $0<\epsilon<1/2$, there exists $\delta=\delta(\epsilon)>0$ such that for $b\le (1-\epsilon)n/2$, Client has a strategy to create a $K_{n^\delta}$-minor. Technically, the proof there adapted the approach of \cite{HKSS08}, briefly described above.

Here we strengthen the result of \cite{HKT16} by proving:
 \begin{thm}\label{thm-CW}
For every $\epsilon>0$ there exists $c>0$ such that for all sufficiently large $n$, when playing a $b$-biased Client--Waiter game on $E(K_n)$ with $b\le (1-\epsilon)\frac{n}{2}$, Client has a strategy to create a minor of the complete graph $K_{c\sqrt{n}}$.
\end{thm}

\noindent{\bf Proof.} Due to bias monotonicity we can assume that $\epsilon$ is small enough where necessary.

Client aims to create a locally sparse bounded degree subgraph $C_1$ on nearly $n$ vertices within his edges; then applying Corollary \ref{complete_minors} allows to argue that $C_1$ contains a minor of $K_{c\sqrt{n}}$. So Client needs to watch out for locally dense pieces, and also for vertices of high degrees. The following technical tool from \cite{DK16} is perfectly suited for this goal; we adapt its formulation here to fit the present circumstances (playing on the edges of $K_n$).

\begin{lemma}\cite{DK16}\label{lemma-DK}
Let $n,b$ be positive integers. Let ${\cal F}$ be a family of subgraphs of the complete graph $K_n$ on $n$ vertices. Assume $\sum_{H\in{\cal F}} (b+1)^{-|E(H)|}<1/2$. Then, when playing the Client--Waiter game  with bias $b$ on the edges of the complete graph on $n$ vertices, Client has a strategy to claim all edges of a graph $C_0$ on $n$ vertices and at least $\left\lfloor\binom{n}{2}/(b+1)\right\rfloor$ edges, not containing any graph from ${\cal F}$.
\end{lemma}

Assume $\delta>0$ is a small constant, whose value will be set soon, and define the two families ${\cal F}_1, {\cal F}_2$  as before. We have:
$$
\sum_{H\in{\cal F}_1} \left(\frac{1}{b+1}\right)^{|E(H)|} \le
\binom{n}{\delta n}\binom{\delta n^2}{\frac{\epsilon n}{4}}\left((1-\epsilon)\frac{n}{2}\right)^{-\frac{\epsilon n}{4}}
\le \left[\left(\frac{e}{\delta}\right)^{\delta}\cdot\left(\frac{30\delta}{\epsilon}\right)^{\frac{\epsilon}{4}}\right]^n= o(1)\,,
$$
for $\delta\le \delta_1(\epsilon)$, and small enough $\epsilon$. Getting to ${\cal F}_2$, we calculate:
$$
\sum_{H\in{\cal F}_2} \left(\frac{1}{b+1}\right)^{|E(H)|} \le
\sum_{k\le \delta n}\binom{n}{k}\binom{\binom{k}{2}}{\left(1+\frac{\epsilon}{8}\right)k}
\left((1-\epsilon)\frac{n}{2}\right)^{-\left(1+\frac{\epsilon}{8}\right)k}= o(1)\,,
$$
for $\delta\le \delta_2$.

 Denoting now ${\cal F}={\cal F}_1\cup{\cal F}_2$, taking $\delta=\min\{\delta_1,\delta_2\}$, and applying Lemma \ref{lemma-DK}, we infer that Client has a strategy to create a graph $C_0$ of his edges, satisfying properties {\bf (P1)}--{\bf (P3)} above. Arguing as in the proof of Theorem \ref{thm-MB}, we conclude that Client gets a graph $C_1$ containing a complete minor $K_t$ with $t=\Theta(\sqrt{n})$. \hfill$\Box$

\bigskip

\noindent{\bf Acknowledgement.} The author wishes to thank Noga Alon and Avi Widgerson for very helpful discussions about the algorithmic aspects of the main result, Dan Hefetz for his input on positional games and his remarks, Wojciech Samotij for his careful reading of an earlier version of the paper.

\end{document}